\documentclass[11pt]{amsart}

\usepackage{amsmath,amssymb,amsthm}
\usepackage{mathtools}
\usepackage[margin=1in]{geometry}
\usepackage[hidelinks]{hyperref}
\usepackage{enumitem}
\usepackage{microtype}
\usepackage{graphicx}
\graphicspath{{figures/}}

\newtheorem{assumption}{Assumption}
\newtheorem{definition}{Definition}
\newtheorem{remark}{Remark}
\newtheorem{lemma}{Lemma}
\newtheorem{theorem}{Theorem}
\newtheorem{corollary}{Corollary}
\newtheorem{fact}{Fact}

\title[LIF Dynamics on Arbitrary Time Scales]{Structure-Preserving Leaky Integrate-and-Fire Dynamics on Arbitrary Time Scales}
\author[B. Kaymak\c{c}alan]{Billur Kaymak\c{c}alan}
\address{University of Turkish Aeronautical Association, Ankara 06790, T\"urkiye}
\email{bkaymakcalan@thk.edu.tr}
\author[B. \.{I}. Sel\c{c}uk]{Behiye \.{I}layda Sel\c{c}uk}
\address{Bilkent University, Ankara 06800, T\"urkiye}
\email{ilayda.selcuk@bilkent.edu.tr}
\date{}

\begin{document}

\begin{abstract}
The leaky integrate-and-fire (LIF) neuron combines a passive membrane
law with threshold-triggered firing and reset. We formulate this
hybrid neuron model, under constant input, on an arbitrary Hilger
time scale $\mathbb{T}$. Our central result is that, after the first
spike, the firing times satisfy
\[
 t_{n+1}=\Pi_{\mathbb{T}}^{+}(t_n+T_c),
\]
where $T_c$ is the classical reset-to-threshold travel time. This
recursion requires no periodicity, homogeneity, or translation
invariance of $\mathbb{T}$; it holds for every time scale that is
unbounded above. To preserve the passive membrane physics, we use
the inverse cylinder transformation to construct a positively
regressive coefficient $p_m$. The resulting delta-dynamic equation
is exactly the restriction to $\mathbb{T}$ of the classical
continuous-time LIF trajectory. By contrast, the direct transcription
$V^\Delta=-(V-V_\infty)/\tau_m$ exhibits graininess-dependent sign
reversal and amplification at coarse right-scattered points. The
first firing time is the forward projection of the classical
continuous threshold-crossing time onto $\mathbb{T}$, and, when the
graininess is bounded by $\mu^\ast$, the per-event forward-projection
delay --- from each reset to the next realized firing time, relative
to the continuous crossing time computed from that same reset --- is
bounded above by $\mu^\ast$. On $\mathbb{T}=h\mathbb{Z}$ we obtain the explicit
lattice-quantized law
$T_{\mathrm{ISI}}^{(h)}=h\lceil T_c/h\rceil$, which converges to the
classical continuous-time inter-spike interval as $h\to0^+$. The
framework therefore separates preserved membrane physics from the
temporal opportunity structure imposed by the time domain. To our
knowledge, this is the first systematic Hilger-calculus treatment of
a threshold-and-reset LIF spiking neuron on an arbitrary time scale.
\end{abstract}

\maketitle

\section{Introduction}

The leaky integrate-and-fire neuron occupies a special position in
computational neuroscience: it is simple enough to admit explicit
analysis, but it retains the essential hybrid structure of a spiking
neuron --- continuous subthreshold voltage evolution together with
threshold-triggered firing and reset. In its classical form, the
model is posed on $\mathbb{R}$ and the passive membrane dynamics are
governed by a first-order linear ordinary differential equation. In
simulations, however, one often works on a regular clock or with
event-driven updates, and in experimental settings the times at
which a state is observed may be nonuniform. This raises a
structural question that is logically prior to any particular
simulation scheme: \emph{can the same membrane law be formulated on
continuous, discrete, and irregular temporal domains without
changing the underlying passive biophysics?}

Hilger's calculus on time scales was introduced precisely to unify
continuous and discrete analysis within a single mathematical
formalism \cite{Hilger1990,Hilger1997,LakshmikanthamSivasundaramKaymakcalan1996,BohnerPeterson2001}. A
\emph{time scale} $\mathbb{T}$ is an arbitrary nonempty closed subset
of $\mathbb{R}$; its local temporal geometry is encoded by the
forward jump operator $\sigma$ and the graininess $\mu$. This makes
time-scale calculus a natural language for the question above.

Time-scale calculus has already been used extensively in
neural-network mathematics, especially for Hopfield
\cite{LiWang2012}, BAM \cite{ChenDu2008,LiGao2010}, Cohen--Grossberg
\cite{LiHuaFei2009}, and other recurrent rate models with delays,
impulses, periodicity, and stability questions. These works concern
network-level state equations rather than the threshold-and-reset
membrane dynamics of a single spiking neuron. A closer neighbor is the work
of Duque, Leiva, Gallo, and Tridane \cite{Duque2023}, whose Example 3
represents a single self-excited neuron (in the sense of Pasemann
\cite{Pasemann1993}) by a semilinear dynamic equation on time scales
and studies bounded, periodic, and almost-periodic solutions. That
model has no firing threshold, spike-generation mechanism, or reset
map, so it does not settle the question we address here; it does,
however, mean that a claim of ``the first single-neuron model on
time scales'' would be incorrect, and we do not make it.

There is also a mature computational-neuroscience literature on
exact or event-driven integration of integrate-and-fire models.
Rotter and Diesmann \cite{RotterDiesmann1999} used matrix
exponentials for exact digital propagation of linear neuronal
systems on a regular grid, and Brette \cite{Brette2006,Brette2007}
developed exact, event-driven integration strategies for broader
classes of integrate-and-fire models. The exponential propagation
formula for a linear LIF membrane is therefore \emph{not itself
new}. Our purpose is different: we place the subthreshold flow,
threshold crossing, and reset into a single Hilger-calculus
framework on an arbitrary $\mathbb{T}$, and derive firing-time and
inter-spike-interval results directly from the geometry of the time
domain, rather than from a simulation scheme.

A second nearby literature concerns state-dependent impulsive neuron
models in classical continuous time. For example, FitzHugh--Nagumo
systems have been augmented with threshold-triggered, state-dependent
reset mechanisms to reproduce spiking and bursting \cite{He2020}.
Such models motivate the hybrid event structure used here, but are
not formulated as delta-dynamic systems on arbitrary Hilger time
scales.

Alternative generalized calculi have also been applied to LIF
dynamics. Fractional-order LIF models, for instance, introduce
nonlocal memory to represent power-law spike-frequency adaptation
\cite{Teka2014}. Our motivation is different: time-scale calculus
does not introduce long-range temporal memory; rather, it unifies
distinct temporal domains and makes their local graininess explicit.

The main contributions of this paper are:
\begin{enumerate}[label=(\roman*)]
\item a structure-preserving Hilger lift of the passive LIF membrane
law, obtained from the classical decay rate by the inverse cylinder
transformation;
\item a complete local classification of the graininess-dependent
failure modes of the direct, or naive, delta-LIF transcription;
\item a forward-projection theorem for threshold-triggered firing on
an arbitrary time scale;
\item a recursive spike-train law valid on every closed time scale
unbounded above, without any periodicity or homogeneity assumption;
\item a graininess-controlled spike-timing bound and an explicit
inter-spike-interval/firing-rate quantization law on $\mathbb{T}=h\mathbb{Z}$.
\end{enumerate}

A distinction that recurs throughout is between \emph{preserved
membrane physics} and \emph{temporal opportunity structure}. In the
structure-preserving formulation, the subthreshold trajectory is
exactly the restriction of the classical continuous solution to
$\mathbb{T}$: the passive biophysics does not change. What changes
with $\mathbb{T}$ is the set of times at which a threshold-triggered
event is permitted to occur. To our knowledge, no previous study has formulated a
threshold-and-reset leaky integrate-and-fire spiking neuron as a
dynamic system on an arbitrary Hilger time scale. We deliberately
make this narrower claim rather than any claim about the first use of
time-scale calculus in neuroscience or the first single-neuron
equation on time scales.

\section{Preliminaries on time scales}\label{sec:prelim}

We recall only the notions used below; see Lakshmikantham, Sivasundaram,
and Kaymak\c{c}alan \cite{LakshmikanthamSivasundaramKaymakcalan1996}
and Bohner and Peterson \cite{BohnerPeterson2001} for the general theory.

A \emph{time scale} $\mathbb{T}$ is a nonempty closed subset of
$\mathbb{R}$. The forward and backward jump operators are
\[
\sigma(t) := \inf\{s \in \mathbb{T} : s > t\}, \qquad
\rho(t) := \sup\{s \in \mathbb{T} : s < t\},
\]
with $\inf \emptyset := \sup \mathbb{T}$ and
$\sup \emptyset := \inf \mathbb{T}$. The \emph{graininess} function is
$\mu(t) := \sigma(t) - t \ge 0$. A point $t$ is right-scattered if
$\sigma(t) > t$ and right-dense if $\sigma(t) = t$.

A function $f : \mathbb{T} \to \mathbb{R}$ is \emph{rd-continuous}
($f \in C_{\mathrm{rd}}$) if it is continuous at right-dense points
and has finite left-sided limits at left-dense points. The
$\Delta$-derivative $f^\Delta(t)$ is defined in the usual way
(Bohner--Peterson, Def.\ 1.10); for right-scattered $t$,
$f^\Delta(t) = \frac{f(\sigma(t)) - f(t)}{\mu(t)}$.

\begin{fact}[Graininess is rd-continuous]\label{fact:mu-rd}
For every time scale $\mathbb{T}$, $\mu \in C_{\mathrm{rd}}(\mathbb{T},[0,\infty))$.
\end{fact}

A function $p \in C_{\mathrm{rd}}(\mathbb{T},\mathbb{R})$ is
\emph{regressive} ($p \in \mathcal{R}$) if $1 + \mu(t)p(t) \neq 0$ for
all $t \in \mathbb{T}$, and \emph{positively regressive}
($p \in \mathcal{R}^+$) if $1+\mu(t)p(t) > 0$ for all $t$.

\begin{definition}[Cylinder transformation]\label{def:cylinder}
We only require the real, positively regressive branch of the Hilger
cylinder transformation. For $h>0$ and $z\in\mathbb{R}$ with
$1+hz>0$, define
\[
\xi_h(z):=\frac{1}{h}\ln(1+hz),
\qquad \xi_0(z):=z.
\]
On this branch $\xi_h$ is a bijection onto $\mathbb{R}$, with inverse
\[
\xi_h^{-1}(w)=\frac{e^{hw}-1}{h},
\qquad \xi_0^{-1}(w)=w.
\]
\end{definition}

\begin{definition}[Generalized exponential]\label{def:exp}
For $p\in\mathcal{R}^+$, the generalized exponential function is
\[
e_p(t,s):=\exp\!\left(
\int_s^t \xi_{\mu(\tau)}(p(\tau))\,\Delta\tau
\right),
\qquad t,s\in\mathbb{T}.
\]
\end{definition}

We use the following standard facts without proof (see Bohner
and Peterson \cite{BohnerPeterson2001}).

\begin{fact}[Cocycle / jump identity]\label{fact:cocycle}
For $p \in \mathcal{R}^+$,
$e_p(\sigma(t),s)=(1+\mu(t)p(t))\,e_p(t,s)$ and
$e_p(t,s)>0$ for all $t\ge s$.
\end{fact}

\begin{fact}[Linear IVP]\label{fact:ivp}
If $p \in \mathcal{R} \cap C_{\mathrm{rd}}(\mathbb{T},\mathbb{R})$, the
initial value problem $y^\Delta(t) = p(t) y(t)$, $y(t_0) = y_0$, has
the unique solution $y(t) = y_0 \, e_p(t,t_0)$ on $\mathbb{T}$.
\end{fact}

\begin{fact}[Integral of the constant function $1$]\label{fact:h1}
For all $s,t \in \mathbb{T}$, $\displaystyle\int_s^t 1 \, \Delta\tau = t - s$.
\end{fact}

Finally, we record the forward projection operator that will encode
firing on $\mathbb{T}$.

\begin{definition}[Forward projection]\label{def:projection}
For $x \in \mathbb{R}$, $\Pi_{\mathbb{T}}^{+}(x) := \inf\{t \in \mathbb{T} : t \ge x\}$.
\end{definition}

\begin{fact}[Well-posedness of $\Pi_{\mathbb{T}}^{+}$]\label{fact:projection-wd}
If $\mathbb{T}$ is unbounded above, then $\Pi_{\mathbb{T}}^{+}(x) \in
\mathbb{T}$ for every $x \in \mathbb{R}$.
\end{fact}
\begin{proof}
$\mathbb{T} \cap [x,\infty)$ is closed (intersection of two closed
sets), nonempty (as $\mathbb{T}$ is unbounded above), and bounded
below by $x$; a nonempty closed subset of $\mathbb{R}$ bounded below
attains its infimum.
\end{proof}

\section{The classical leaky integrate-and-fire neuron}\label{sec:classical}

The subthreshold membrane dynamics of the LIF neuron under constant
input $I_0$ is
\[
C_m \dot V(t) = -g_L\big(V(t) - E_L\big) + I_0 .
\]
Here $V$ is the membrane potential, $C_m$ the membrane capacitance,
$g_L$ the leak conductance, $E_L$ the leak (resting) potential, and
$I_0$ a constant injected current. Writing $\tau_m := C_m/g_L$ and
$V_\infty := E_L + I_0/g_L$,
\begin{equation}\label{eq:classical-lif}
\dot V(t) = -\frac{1}{\tau_m}\big(V(t) - V_\infty\big).
\end{equation}
The parameter $V_{th}$ denotes the firing threshold and $V_r$ the
post-spike reset potential. A spike is declared when the current
unreset branch first reaches $V_{th}$, after which the state is reset
to $V_r<V_{th}$. Given $V_\infty > V_{th} > V(t_0) =: V_0$, the
classical crossing time from $\eqref{eq:classical-lif}$ is
\begin{equation}\label{eq:classical-crossing}
t^\ast = t_0 + \tau_m \ln\!\left(\frac{V_\infty - V_0}{V_\infty - V_{th}}\right),
\end{equation}
and, given $V_\infty > V_{th} > V_r$, the reset-to-threshold travel
time is the constant
\begin{equation}\label{eq:Tc}
T_c := \tau_m \ln\!\left(\frac{V_\infty - V_r}{V_\infty - V_{th}}\right).
\end{equation}

\begin{assumption}[Standing hypotheses]\label{ass:standing}
Throughout the paper, unless stated otherwise:
\begin{enumerate}[label=(A\arabic*)]
\item\label{A1} $\mathbb{T}$ is a time scale, unbounded above, and
$t_0 \in \mathbb{T}$.
\item\label{A2} $\tau_m > 0$ is a constant.
\item\label{A3} The input is constant, so $V_\infty$ is a constant.
\item\label{A4} (Sections~\ref{sec:firing}--\ref{sec:isi} only)
$V_\infty > V_{th} > V_r$, and $V_0 < V_{th}$.
\item\label{A5} (Section~\ref{sec:isi} only)
$\mu^\ast := \sup_{t \in \mathbb{T}} \mu(t) < \infty$.
\end{enumerate}
Beyond \ref{A1}, $\mathbb{T}$ is \emph{not} assumed homogeneous,
periodic, or translation-invariant anywhere in this paper.
\end{assumption}

\section{Two time-scale extensions of the subthreshold equation}\label{sec:extensions}

\subsection{The naive delta-LIF equation}

\begin{definition}[Naive delta-LIF]\label{def:naive}
\[
V^\Delta(t) = -\frac{1}{\tau_m}\big(V(t) - V_\infty\big).
\]
\end{definition}

At a right-scattered point, Definition~\ref{def:naive} gives
\begin{equation}\label{eq:naive-step}
V(\sigma(t)) - V_\infty = M(t)\,\big(V(t) - V_\infty\big), \qquad
M(t) := 1 - \frac{\mu(t)}{\tau_m}.
\end{equation}

\begin{theorem}[Local failure regimes of the naive $\Delta$-LIF equation]\label{thm:1}
Let $t \in \mathbb{T}$ be right-scattered. With $M(t)$ as in
\eqref{eq:naive-step}:
\begin{enumerate}[label=(\roman*)]
\item If $0 < \mu(t) < \tau_m$: $0 < M(t) < 1$ (sign-preserving decay).
\item If $\mu(t) = \tau_m$: $M(t) = 0$; regressivity of the naive
coefficient $-1/\tau_m$ fails at $t$ ($1+\mu(t)(-1/\tau_m)=0$).
\item If $\tau_m < \mu(t) < 2\tau_m$: $-1 < M(t) < 0$ (sign-alternating,
decaying magnitude).
\item If $\mu(t) = 2\tau_m$: $M(t) = -1$ (neutral, undamped alternation).
\item If $\mu(t) > 2\tau_m$: $|M(t)| > 1$ (sign-alternating, growing
magnitude; local amplification).
\end{enumerate}
Consequently, if $\mathbb{T}=h\mathbb{Z}$ with $h>2\tau_m$ and
$V_0\neq V_\infty$, then
$|V(t_0+nh)-V_\infty|=|M|^n|V_0-V_\infty|\to\infty$. (If
$V_0=V_\infty$ the naive trajectory is identically $V_\infty$, by
\eqref{eq:naive-step}, and there is nothing to diverge.) More
generally, suppose there is a sequence of consecutive right-scattered
points $t_{k+1}=\sigma(t_k)$ with $t_k\to\infty$,
$\inf_k\mu(t_k)>2\tau_m$, and $V(t_0)\neq V_\infty$ at the first
point of that sequence. Then $|V(t_k)-V_\infty|$ grows at least
geometrically along that sequence and is therefore unbounded.
\end{theorem}
\begin{proof}
(i)--(v) are immediate from \eqref{eq:naive-step}. For the
$h\mathbb{Z}$ corollary, $M(t) \equiv 1-h/\tau_m$ is constant, so the
$n$-step propagator is $M^n$. For the general statement, set
$q:=\inf_k(\mu(t_k)/\tau_m-1)$. The hypothesis
$\inf_k\mu(t_k)>2\tau_m$ gives $q>1$, while
$|M(t_k)|=\mu(t_k)/\tau_m-1\ge q$. Hence
\[
|V(t_n)-V_\infty|\ge q^n|V(t_0)-V_\infty|,
\]
along the consecutive scattered sequence (after re-indexing its first
point as $t_0$). Since $V(t_0)\neq V_\infty$ by hypothesis, the
right-hand side is a nonzero constant times $q^n$ with $q>1$, and
therefore diverges.
\end{proof}

\subsection{The structure-preserving LIF equation}

\begin{definition}[Structure-preserving LIF]\label{def:sp}
\[
p_m(t) := \xi_{\mu(t)}^{-1}\!\left(-\frac{1}{\tau_m}\right)
= \begin{cases}
\dfrac{e^{-\mu(t)/\tau_m} - 1}{\mu(t)}, & \mu(t) > 0, \\[1.2ex]
-\dfrac{1}{\tau_m}, & \mu(t) = 0,
\end{cases}
\qquad
V^\Delta(t) = p_m(t)\big(V(t) - V_\infty\big).
\]
We emphasize that $p_m$ is the direct application of the inverse
cylinder transformation (Definition~\ref{def:cylinder}) to the
constant continuous-time rate $-1/\tau_m$; it is not itself a new
piece of time-scale machinery.
\end{definition}

\begin{remark}[$p_m$ is rd-continuous]\label{rem:1}
By Fact~\ref{fact:mu-rd}, $\mu \in C_{\mathrm{rd}}(\mathbb{T},[0,\infty))$.
The map $g(x) := (e^{-x/\tau_m}-1)/x$ for $x>0$, $g(0):=-1/\tau_m$, is
continuous on $[0,\infty)$ (removable singularity at $0$, matching
the piecewise definition above). Hence $p_m = g \circ \mu \in
C_{\mathrm{rd}}(\mathbb{T},\mathbb{R})$, with no further hypothesis
on $\mathbb{T}$.
\end{remark}

\begin{lemma}[Positive regressivity]\label{lem:1}
$1 + \mu(t) p_m(t) = e^{-\mu(t)/\tau_m} > 0$ for all $t \in \mathbb{T}$;
in particular $p_m \in \mathcal{R}^+$.
\end{lemma}
\begin{proof}
For $\mu(t)>0$: $1 + \mu(t)\cdot\frac{e^{-\mu(t)/\tau_m}-1}{\mu(t)} =
e^{-\mu(t)/\tau_m}$. For $\mu(t)=0$: $1+0 = 1 = e^0$.
\end{proof}

\begin{lemma}[Hilger exponential identity]\label{lem:2}
$\xi_{\mu(t)}(p_m(t)) \equiv -1/\tau_m$ for all $t \in \mathbb{T}$,
and consequently $e_{p_m}(t,s) = e^{-(t-s)/\tau_m}$ for all $t,s \in
\mathbb{T}$.
\end{lemma}
\begin{proof}
The first identity holds by construction ($p_m(t) =
\xi_{\mu(t)}^{-1}(-1/\tau_m)$, so applying $\xi_{\mu(t)}$ recovers
$-1/\tau_m$). By Definition~\ref{def:exp} and Fact~\ref{fact:h1},
\[
e_{p_m}(t,s) = \exp\left(\int_s^t \left(-\frac{1}{\tau_m}\right)\Delta\tau\right)
= \exp\left(-\frac{t-s}{\tau_m}\right). \qedhere
\]
\end{proof}

\begin{corollary}[Exact restriction of the classical trajectory]\label{cor:1}
Under \ref{A1}--\ref{A3}, the unique solution of $V^\Delta(t) =
p_m(t)(V(t)-V_\infty)$, $V(t_0) = V_0$, is
\begin{equation}\label{eq:restriction}
V(t) = V_\infty + (V_0 - V_\infty)\, e^{-(t-t_0)/\tau_m}, \qquad t \in \mathbb{T},\ t \ge t_0.
\end{equation}
That is, $V|_{\mathbb{T}}$ is the exact restriction to $\mathbb{T}$
of the classical continuous-time solution of \eqref{eq:classical-lif}.
\end{corollary}
\begin{proof}
Let $y := V - V_\infty$ (so $y^\Delta = V^\Delta$ since $V_\infty$ is
constant); then $y^\Delta(t) = p_m(t)y(t)$, $y(t_0)=V_0-V_\infty$. By
Lemma~\ref{lem:1} and Remark~\ref{rem:1}, $p_m \in \mathcal{R}^+ \cap
C_{\mathrm{rd}}$, so Fact~\ref{fact:ivp} applies:
$y(t) = (V_0-V_\infty)\,e_{p_m}(t,t_0)$. Apply Lemma~\ref{lem:2}.
\end{proof}

\begin{corollary}[Monotonicity; no artificial overshoot]\label{cor:2}
Under Corollary~\ref{cor:1}, if $V_0 < V_\infty$ then $V_0 \le V(t) <
V_\infty$ for all $t \in \mathbb{T}$, $t \ge t_0$; symmetrically, if
$V_0 > V_\infty$ then $V_\infty < V(t) \le V_0$. This holds for
\emph{every} $\mathbb{T}$ and every graininess profile -- no
condition on $\mu$ is required.
\end{corollary}
\begin{proof}
For $t \ge t_0$, $e^{-(t-t_0)/\tau_m} \in (0,1]$ regardless of
$\mathbb{T}$. If $V_0 < V_\infty$ then $V_0 - V_\infty < 0$, so
$V(t)-V_\infty = (V_0-V_\infty)e^{-(t-t_0)/\tau_m} \in [V_0-V_\infty,\,0)$
by \eqref{eq:restriction}, i.e.\ $V(t) \in [V_0, V_\infty)$. The
symmetric case is identical.
\end{proof}

\section{Firing and reset on an arbitrary time scale}\label{sec:firing}

Assume \ref{A4}. By \eqref{eq:classical-crossing}, the classical
continuous solution $V_c(t) = V_\infty - (V_\infty-V_0)e^{-(t-t_0)/\tau_m}$
is strictly increasing (as $V_\infty > V_0$), converges to $V_\infty$,
and satisfies $V_c(t) < V_{th}$ for $t<t^\ast$, $V_c(t) \ge V_{th}$
for $t \ge t^\ast$.

\begin{theorem}[First firing time as forward projection]\label{thm:2}
Under \ref{A1} and \ref{A4}, the first time at which the
structure-preserving trajectory \eqref{eq:restriction} satisfies
$V(t) \ge V_{th}$ is
\[
t_f = \Pi_{\mathbb{T}}^{+}(t^\ast).
\]
\end{theorem}
\begin{proof}
By Corollary~\ref{cor:1}, $V(t) = V_c(t)$ for all $t \in \mathbb{T}$
with $t\ge t_0$ (the restriction property places no threshold
restriction on its validity as an identity of real numbers; the
subthreshold label is purely interpretive). Hence, for $t \in
\mathbb{T}$: $V(t) < V_{th} \iff t < t^\ast$ and $V(t) \ge V_{th}
\iff t \ge t^\ast$. The first such $t$ is, by definition,
$\Pi_{\mathbb{T}}^{+}(t^\ast)$, which lies in $\mathbb{T}$ by
Fact~\ref{fact:projection-wd} and \ref{A1}.
\end{proof}

At a firing time we distinguish the value of the current
\emph{unreset branch}, which is used to detect threshold attainment,
from the post-event state used to initialize the next branch. Thus,
at $t_f$ the unreset branch satisfies $V_{\mathrm{pre}}(t_f)\ge V_{th}$,
a spike is recorded, and the post-event state is set to
$V_{\mathrm{post}}(t_f)=V_r$. The next subthreshold branch is then
initialized at the same time $t_f$ with value $V_r$. These
``pre''/``post'' labels refer to event-side states, not to the
forward/backward jump operators of the time scale. Write $t_1:=t_f$
and let $t_2,t_3,\dots$ denote the subsequent firing times.

\begin{theorem}[Recursive spike train on an arbitrary time scale]\label{thm:3}
Under \ref{A1} and \ref{A4}, the firing times obey
\[
t_{n+1} = \Pi_{\mathbb{T}}^{+}(t_n + T_c), \qquad n = 1,2,\dots,
\]
with $T_c$ as in \eqref{eq:Tc}. This holds for \emph{every} time
scale satisfying \ref{A1} -- no periodicity, homogeneity, or
translation-invariance of $\mathbb{T}$ is required.
\end{theorem}
\begin{proof}
Let $t_1:=t_f$, whose existence is given by Theorem~\ref{thm:2}.
Fix any $n\ge1$ for which the firing time $t_n$ has been generated.
Immediately after the event, $V_{\mathrm{post}}(t_n)=V_r$; we denote
the restarted subthreshold branch again by $V$, with initial
condition $V(t_n)=V_r$. Because \eqref{eq:classical-lif} is
autonomous ($\tau_m$ and $V_\infty$ are constant),
Corollary~\ref{cor:1}, restarted at $t_n$, gives
\[
V(t)=V_\infty+(V_r-V_\infty)e^{-(t-t_n)/\tau_m},
\qquad t\in\mathbb{T},\ t\ge t_n.
\]
Hence the corresponding continuous reset-to-threshold crossing time
is $t_n+T_c$, independent of the absolute position of $t_n$ in
$\mathbb{T}$. Applying the forward-projection argument of
Theorem~\ref{thm:2} to this restarted branch yields
\[
t_{n+1}=\Pi_{\mathbb{T}}^+(t_n+T_c),
\]
which is well posed by Fact~\ref{fact:projection-wd}. Repeating this
construction after each reset generates the stated recursion for all
$n\ge1$.
\end{proof}

\begin{remark}
Theorem~\ref{thm:3} is the central result of the paper: what makes the reset-to-threshold travel time independent of
absolute time is the autonomy of the subthreshold ODE, not any
regularity of $\mathbb{T}$. The time scale
may be continuous, uniformly discrete, or arbitrarily irregular --
the same one-line recursion governs the spike train in every case.
In particular, $\Pi_{\mathbb{T}}^{+}(x)\ge x$ for every $x$, so
$t_{n+1}\ge t_n+T_c$ and hence $t_n\ge t_1+(n-1)T_c\to\infty$; the
model therefore admits no finite-time accumulation of firing events
(no Zeno behavior).
\end{remark}

\section{Spike-timing error and firing-rate quantization}\label{sec:isi}

\begin{lemma}[General projection-gap bound]\label{lem:3}
Let $\mathbb{T}$ satisfy \ref{A1} and \ref{A5}, and let $x \in
\mathbb{R}$ with $\mathbb{T} \cap (-\infty,x] \neq \emptyset$.
Define $\delta_{\mathbb{T}}(x) := \Pi_{\mathbb{T}}^{+}(x) - x$. Then
\[
0 \le \delta_{\mathbb{T}}(x) \le \mu^\ast .
\]
This is a statement about $\mathbb{T}$ alone; it does not reference
the neuron model.
\end{lemma}
\begin{proof}
$\delta_{\mathbb{T}}(x) \ge 0$ is immediate from
Definition~\ref{def:projection}. If $x \in \mathbb{T}$, $\delta=0$.
Otherwise let $t := \sup(\mathbb{T} \cap (-\infty,x])$; since
$\mathbb{T}$ is closed and this set is nonempty and bounded above by
$x$, $t \in \mathbb{T}$ and $t \le x$. Let $s :=
\Pi_{\mathbb{T}}^{+}(x)$. By maximality of $t$ and minimality of $s$,
$(t,s) \cap \mathbb{T} = \emptyset$, so $s = \sigma(t)$ and
$\mu(t) = s-t$. Since $t \le x < s$,
$\delta_{\mathbb{T}}(x) = s-x \le s-t = \mu(t) \le \mu^\ast$.
\end{proof}

\begin{corollary}[Spike-time error bound]\label{cor:3}
Under \ref{A1}, \ref{A4}, \ref{A5}: $0 \le t_f - t^\ast \le \mu^\ast$,
and more generally $0 \le t_{n+1} - (t_n+T_c) \le \mu^\ast$ for all
$n$.
\end{corollary}
\begin{proof}
Apply Lemma~\ref{lem:3} with $x=t^\ast$ (resp.\ $x=t_n+T_c$).
\end{proof}

\begin{corollary}[$h\mathbb{Z}$ inter-spike interval and firing rate]\label{cor:4}
Let $\mathbb{T} = h\mathbb{Z}$, $h>0$, and suppose a spike occurs at
$t_n \in h\mathbb{Z}$. Then
\[
T_{\mathrm{ISI}}^{(h)} := t_{n+1}-t_n = h\left\lceil \frac{T_c}{h} \right\rceil,
\qquad
f_h := \frac{1}{T_{\mathrm{ISI}}^{(h)}} = \frac{1}{h\lceil T_c/h \rceil}.
\]
As $h \to 0^+$, $T_{\mathrm{ISI}}^{(h)} \to T_c$ and $f_h \to 1/T_c$,
recovering the classical firing rate.
\end{corollary}
\begin{proof}
$\Pi_{h\mathbb{Z}}^{+}(x) = h\lceil x/h \rceil$; apply
Theorem~\ref{thm:3} with $t_n \in h\mathbb{Z}$, giving $t_{n+1} =
h\lceil (t_n+T_c)/h\rceil = t_n + h\lceil T_c/h\rceil$ (as
$t_n/h\in\mathbb{Z}$). The limit is immediate since $h\lceil
T_c/h\rceil - T_c \in [0,h) \to 0$.
\end{proof}

\section{Numerical illustrations}\label{sec:numerics}
\begin{figure}[htbp]
\centering
\includegraphics[width=\textwidth]{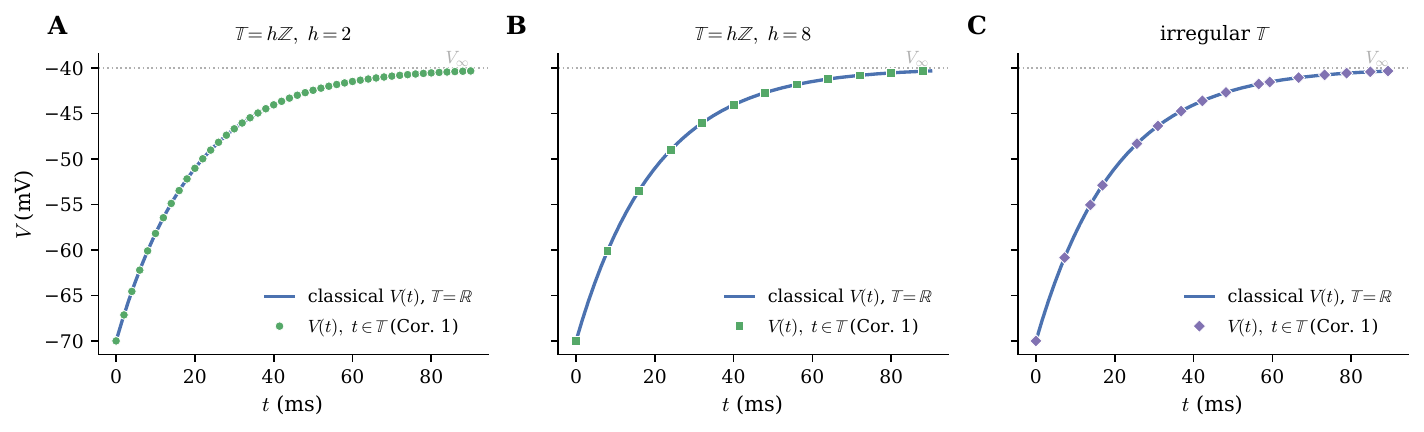}
\caption{Exact restriction of the classical LIF trajectory to
different time scales. The classical continuous-time passive,
unreset trajectory $V(t)$ (drawn beyond threshold, since
Corollary~\ref{cor:1} places no threshold restriction on its
validity) is shown together with the structure-preserving
solution \eqref{eq:restriction} evaluated on three representative
time scales: (A) the uniform lattice $\mathbb{T}=2\mathbb{Z}$~ms, (B)
the coarser lattice $\mathbb{T}=8\mathbb{Z}$~ms, and (C) a
deterministic irregular time scale. In every case the time-scale
solution lies exactly on the classical continuous trajectory at all
$t\in\mathbb{T}$, illustrating Corollary~\ref{cor:1}. Thus the
structure-preserving formulation changes the admissible temporal
domain without altering the passive membrane trajectory itself.}
\label{fig:1}
\end{figure}

\begin{figure}[htbp]
\centering
\includegraphics[width=\textwidth]{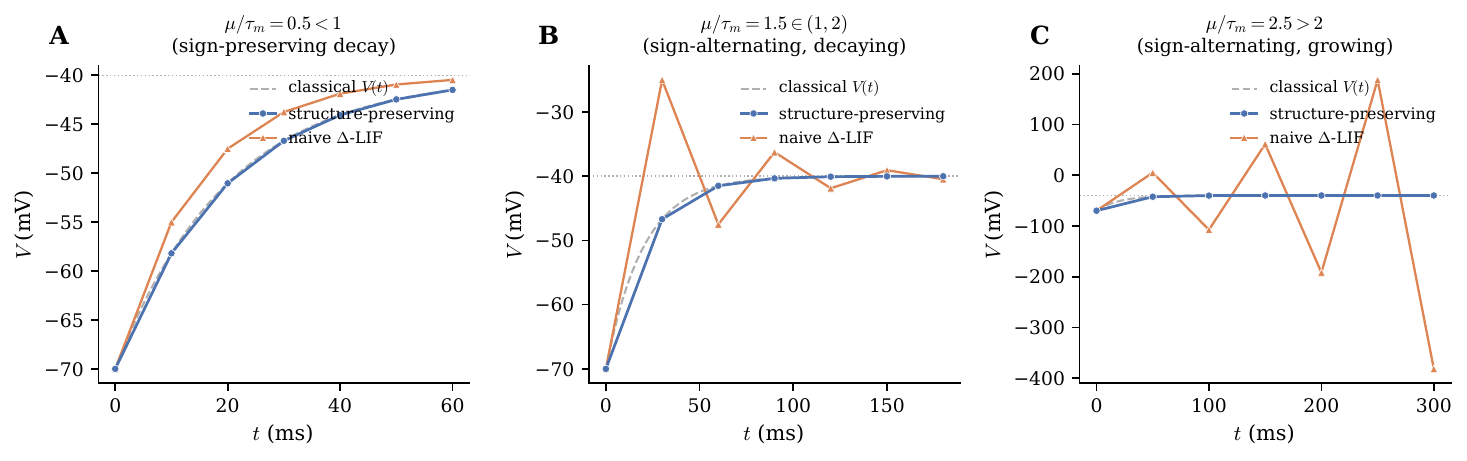}
\caption{Graininess-dependent failure of the naive delta-LIF
formulation. Naive and structure-preserving LIF trajectories are
compared on uniform time scales with three values of the
dimensionless graininess ratio $\rho:=\mu/\tau_m$. (A) For
$\rho=0.5<1$, both formulations decay toward $V_\infty$, although the
naive trajectory does not exactly reproduce the continuous membrane
semigroup. (B) For $\rho=1.5\in(1,2)$, the naive multiplier becomes
negative, producing sign-alternating overshoots around $V_\infty$,
whereas the structure-preserving trajectory remains monotone
(Corollary~\ref{cor:2}). (C) For $\rho=2.5>2$, the magnitude of the
naive multiplier exceeds one and deviations from equilibrium grow,
while the structure-preserving solution remains bounded and
convergent. Threshold/reset events are suppressed in these panels to
isolate the passive dynamics; the panels illustrate the contrast
between Theorem~\ref{thm:1} and Corollary~\ref{cor:2}.}
\label{fig:2}
\end{figure}

\begin{figure}[htbp]
\centering
\includegraphics[width=0.85\textwidth]{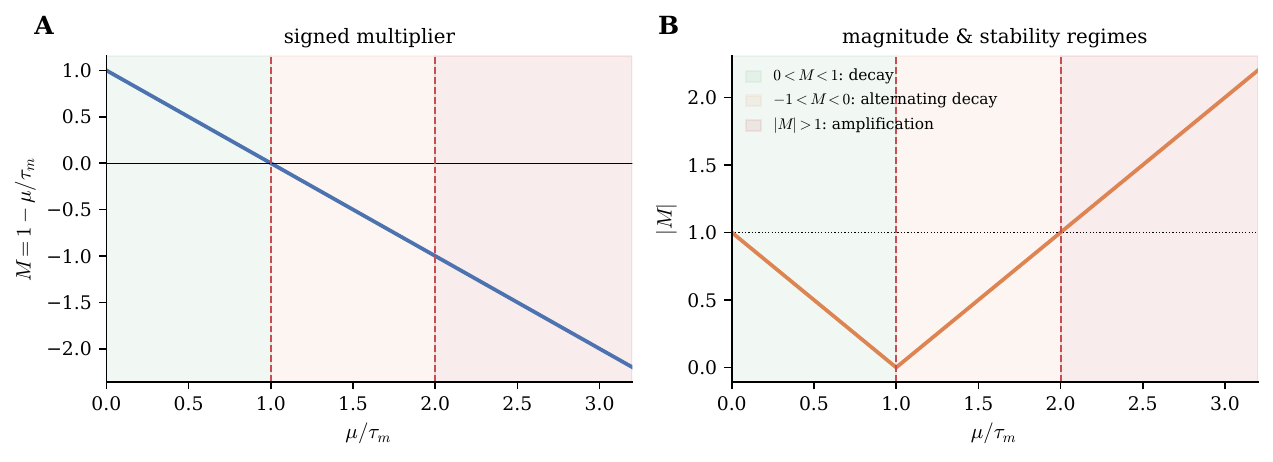}
\caption{Local dynamical regimes of the naive delta-LIF multiplier.
The naive right-scattered update \eqref{eq:naive-step} is governed by
the dimensionless multiplier $M(\rho)=1-\rho$, $\rho:=\mu/\tau_m$.
(A) Signed multiplier $M(\rho)$. (B) Magnitude $|M(\rho)|$. For
$0<\rho<1$ the update is sign-preserving and contractive; for
$1<\rho<2$ it is sign-alternating but contractive; at $\rho=2$,
$|M|=1$; and for $\rho>2$, $|M|>1$, resulting in local amplification.
The vertical boundaries at $\rho=1$ and $\rho=2$ correspond to the
regime transitions characterized in Theorem~\ref{thm:1}.}
\label{fig:3}
\end{figure}

\begin{figure}[htbp]
\centering
\includegraphics[width=0.9\textwidth]{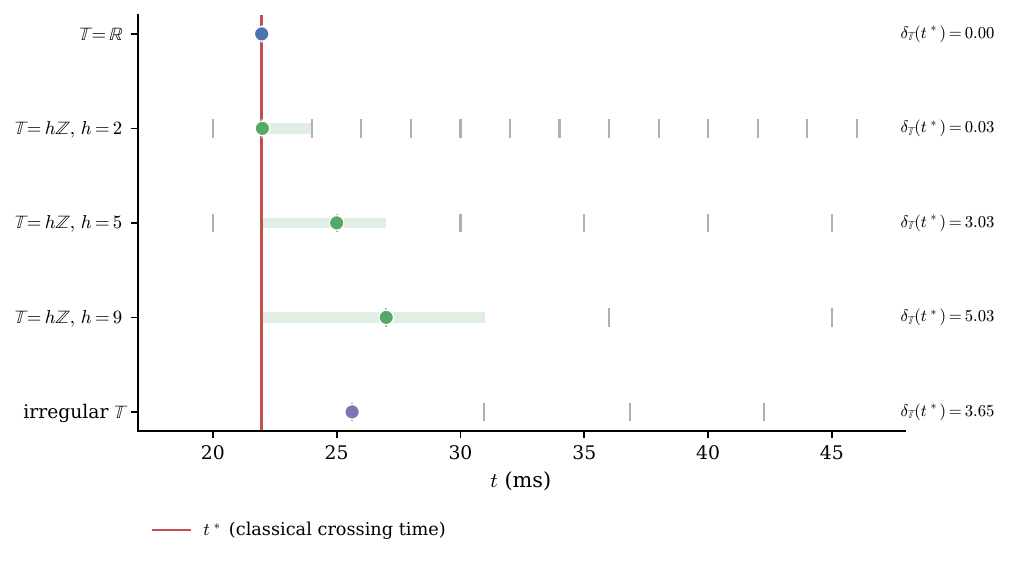}
\caption{Forward projection of the continuous threshold-crossing time
onto a time scale. The vertical line marks the classical
continuous-time threshold-crossing time $t^\ast$. For each temporal
domain, the filled marker denotes the realized on-time-scale firing
time $t_f=\Pi_{\mathbb{T}}^{+}(t^\ast)$ (Theorem~\ref{thm:2}). For
$\mathbb{T}=\mathbb{R}$, the projection is exact and $t_f=t^\ast$. On
uniform lattices $h\mathbb{Z}$, firing is delayed until the first
admissible time-scale point at or after $t^\ast$. The corresponding
projection error $\delta_{\mathbb{T}}(t^\ast)=t_f-t^\ast$ is
annotated for each time scale. For the uniform lattices, the shaded
interval represents the graininess bound $\mu^\ast=h$, illustrating
$0\le\delta_{\mathbb{T}}(t^\ast)\le\mu^\ast$ from
Corollary~\ref{cor:3}. The irregular example illustrates the same
forward-projection mechanism without imposing periodicity or
homogeneity on $\mathbb{T}$. For the deterministic irregular example,
$\mu^\ast$ is also finite and computable (here $\mu^\ast \le 9.3$~ms),
but no graininess band is drawn for visual clarity because a single
uniform interval does not represent its locally varying graininess as
intuitively as for $h\mathbb{Z}$.}
\label{fig:4}
\end{figure}

\begin{figure}[htbp]
\centering
\includegraphics[width=\textwidth]{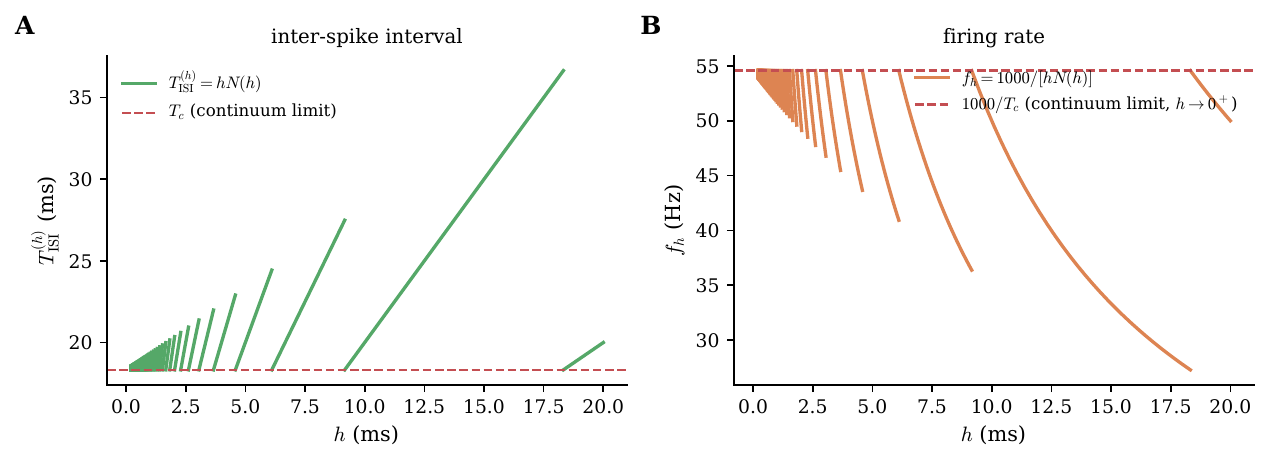}
\caption{Quantization of the inter-spike interval and firing rate on
$h\mathbb{Z}$. The genuine staircase quantity is the integer step
count $N(h):=\lceil T_c/h\rceil$. The plotted quantities
$T_{\mathrm{ISI}}^{(h)}=hN(h)$ and $f_h=1000/(hN(h))$~Hz are
therefore lattice-quantized rather than step functions themselves: on
intervals where $N(h)$ is constant they vary linearly and
hyperbolically in $h$, respectively (with $h$ in ms, for the
representative parameters above), with a jump discontinuity each
time $N(h)$ decreases. (A) $T_{\mathrm{ISI}}^{(h)}$
(Corollary~\ref{cor:4}); the dashed horizontal line indicates the
continuous reset-to-threshold travel time $T_c$. (B) The
corresponding firing rate $f_h$. As $h\to0^+$,
$T_{\mathrm{ISI}}^{(h)}\to T_c$ and $f_h\to1000/T_c$~Hz, recovering
the classical continuous-time firing rate.}
\label{fig:5}
\end{figure}

The numerical examples below are designed only to illustrate the
analytical results above; no parameter is fitted to data and no
theorem depends on the simulations. For concreteness, the following
representative parameter values are used throughout this section:
$\tau_m=20$~ms, $E_L=V_0=-70$~mV, $V_\infty=-40$~mV, $V_{th}=-50$~mV,
$V_r=-65$~mV (so $T_c\approx18.33$~ms by \eqref{eq:Tc}). All five
figures are produced by a single reproducible script
(\texttt{make\_figures.py}); the irregular time scale used in
Figures~\ref{fig:1} and~\ref{fig:4} is a fixed deterministic,
quasi-periodically modulated point sequence, not a random draw. In
Figures~\ref{fig:1}--\ref{fig:3}, threshold/reset events are
intentionally suppressed so that the unreset passive membrane branch
can be displayed beyond threshold; firing semantics are restored in
Figures~\ref{fig:4}--\ref{fig:5}.

\section{Discussion}\label{sec:discussion}

The results above separate two effects that are often conflated in
discrete or irregularly timed neuron simulations. First, the
structure-preserving coefficient $p_m$ makes the passive membrane
trajectory independent of the local graininess in a precise sense:
for every $t \in \mathbb{T}$, the voltage equals the classical
continuous LIF voltage at the same physical time
(Corollary~\ref{cor:1}). Changing $\mathbb{T}$ does not alter the
passive membrane semigroup. Second, the spike train generally
\emph{does} depend on $\mathbb{T}$ under the on-time-scale firing semantics of
Section~\ref{sec:firing}: the continuous
threshold time $t^\ast$ is a latent physical crossing time, but the
realized firing time is its forward projection onto the admissible
temporal domain (Theorem~\ref{thm:2}). The geometry of $\mathbb{T}$
therefore acts as a \emph{temporal opportunity structure} for
events, in contrast to the \emph{preserved membrane physics} of
Corollary~\ref{cor:1}. This distinction is especially transparent on
$h\mathbb{Z}$, where the firing rate becomes the piecewise-hyperbolic,
lattice-quantized function of Corollary~\ref{cor:4} even though the
subthreshold voltage law
itself remains the exact restriction of the continuous solution.

This viewpoint also clarifies the role of the naive $\Delta$-LIF
equation. Equation~\eqref{eq:naive-step} preserves the differential
equation's coefficient but not its continuous semigroup; at coarse
right-scattered points this can produce sign reversal or
amplification (Theorem~\ref{thm:1}) that has no counterpart in the
passive RC membrane. By contrast, the structure-preserving equation
preserves the semigroup exactly, through positive regressivity
(Lemma~\ref{lem:1}). On the uniform lattice $\mathbb{T}=h\mathbb{Z}$,
the naive delta-LIF update is algebraically identical to forward
Euler applied to the passive membrane ODE~\eqref{eq:classical-lif};
thus the $h=2\tau_m$ boundary in Theorem~\ref{thm:1} is the familiar
Euler stability threshold, and Theorem~\ref{thm:1} expresses the same
mechanism locally, in terms of graininess, on an arbitrary time
scale.

\textbf{Relation to exact/event-driven LIF integration.} The
inverse-cylinder coefficient $p_m$ is standard time-scale machinery;
the contribution is not the discovery of the scalar exponential
update, which is well established in computational neuroscience
\cite{RotterDiesmann1999,Brette2006,Brette2007}. What is developed
here is the unified Hilger-calculus treatment of membrane
propagation, threshold attainment, reset, and spike timing on an
arbitrary $\mathbb{T}$ --- a structural, rather than computational,
question.

\textbf{Relation to time-scale neural-network models.} Hopfield
\cite{LiWang2012}, BAM \cite{ChenDu2008,LiGao2010}, Cohen--Grossberg
\cite{LiHuaFei2009}, and related recurrent models on time scales
primarily address network-level existence, periodicity,
synchronization, and stability questions.
Duque et al.\ \cite{Duque2023} give a genuine single-neuron voltage
example on time scales, so a claim of ``the first single-neuron
model on time scales'' would be incorrect; the narrower
threshold-and-reset claim stated in the Introduction is the relevant
one for this work.

\textbf{Relation to state-dependent impulsive neuron models.} A
nearby classical-time literature augments continuous neuron models,
such as FitzHugh--Nagumo, with threshold-triggered, state-dependent
reset mechanisms to reproduce spiking and bursting \cite{He2020}.
Such models motivate the hybrid event structure used here but are
not formulated as delta-dynamic systems on arbitrary Hilger time
scales.

\textbf{On-time-scale versus event-enriched firing semantics.} Under
the on-time-scale firing convention used in Theorem~\ref{thm:2}, with
firing times encoded by the forward projection of
Definition~\ref{def:projection}, a spike can occur only at a point of
$\mathbb{T}$. An alternative,
\emph{event-enriched}, semantics would insert each continuous
threshold-crossing time exactly into the admissible domain (for a
single episode, $\mathbb{T}^\ast=\mathbb{T}\cup\{t^\ast\}$).
Repeated firing would therefore generate a state-dependent enrichment
of the time domain. This construction is closer in spirit to precise
event-driven simulation, but it raises a different class of
mathematical questions; we leave state-dependent time-scale
enrichment for future work.

\textbf{Limitations and future work.} The present formulation is
intentionally minimal. The constant-input LIF model is used because
it permits a complete analytical separation between passive membrane
propagation and the temporal geometry of threshold opportunities.
Time-varying input $I(t)$, refractory intervals, synaptic currents,
adaptive thresholds, nonlinear membrane models (FitzHugh--Nagumo,
Hodgkin--Huxley), and network coupling are natural extensions left
to future work; the constant-input theory developed here should be
viewed as the foundational case in which the role of $\mathbb{T}$
can be characterized exactly. Periodic time scales may also induce
additional arithmetic structure in the inter-spike-interval
sequence; we leave its characterization for future work.

\end{document}